\documentclass{article}

\usepackage{myarxiv}
\usepackage[utf8]{inputenc} 
\usepackage[T1]{fontenc}    
\usepackage{hyperref}       
\usepackage{url}            
\usepackage{booktabs}       
\usepackage{amsfonts}       
\usepackage{nicefrac}       
\usepackage{microtype}      
\usepackage{lipsum}

\usepackage{amssymb,amsmath,amsfonts,amsthm,enumerate}
\usepackage{subfigure}
\usepackage{graphicx}
\usepackage{epsfig}
\usepackage{float}
\usepackage[usenames]{color}
\usepackage{cite}
\usepackage{bm}
\usepackage{enumitem}
\setlist{nosep}

\newtheorem{lemma}{Lemma}[section]

\newtheorem{proposition}{Proposition}[section]
\newtheorem{remark}{Remark}[section]
\newtheorem{theorem}{Theorem}[section]

\def\BC{\mathbb C}
\def\BD{\mathbb D}
\def\BN{\mathbb N}
\def\BR{\mathbb R}

\def\cD{\mathcal D}

\def\cL{\mathcal L}

\def\cU{\mathcal U}

\def\rd{\mathrm d}

\def\rdiv{\mathrm{div}}

\def\Ga{\Gamma}

\def\Om{\Omega}
\def\al{\alpha}
\def\be{\beta}
\def\ga{\gamma}
\def\de{\delta}

\def\la{\lambda}

\def\vp{\varphi}
\def\om{\omega}
\def\f{\frac}
\def\nb{\nabla}
\def\ov{\overline}
\def\pa{\partial}

\def\tri{\triangle}

\title{Inverse Problems of Determining Parameters of the Fractional Partial Differential Equations}

\author{
  Zhiyuan Li\\
  School of Mathematics and Statistics\\
  Shandong University of Technology\\
  Zibo, Shandong 255049, China\\
  \texttt{zyli@sdut.edu.cn}\\
  \And
  Yikan Liu\\
  Graduate School of Mathematical Sciences\\
  The University of Tokyo\\
  3-8-1 Komaba, Meguro-ku, Tokyo 153-8914, Japan\\
  \texttt{ykliu@ms.u-tokyo.ac.jp}\\
  \AND
  Masahiro Yamamoto\\
  Graduate School of Mathematical Sciences\\
  The University of Tokyo\\
  3-8-1 Komaba, Meguro-ku, Tokyo 153-8914, Japan\\
  \texttt{myama@ms.u-tokyo.ac.jp}\\
}

\begin{document}
\maketitle

\begin{abstract}
When considering fractional diffusion equation as model equation in analyzing anomalous diffusion processes, some important parameters in the model related to orders of the fractional derivatives, are often unknown and difficult to be directly measured, which requires one to discuss inverse problems of identifying these physical quantities from some indirectly observed information of solutions. Inverse problems in determining these unknown parameters of the model are not only theoretically interesting, but also necessary for finding solutions to initial-boundary value problems and studying properties of solutions. This chapter surveys works on such inverse problems for fractional diffusion equations.
\end{abstract}

\MRsubject{35R11\and 26A33\and 35R30\and 65M32}

\section{Introduction}\label{sec-intro-multi}

The fractional equations have been playing important roles in various fields such as physics, chemistry, astrophysics during the last few decades. In particular, in heterogenous media, the diffusion often indicates anomalous profiles which cannot be simulated by the classical diffusion equation. Thus several model equations have been introduced, and a time-fractional diffusion equation is one of them. For a flexible type of time-fractional diffusion equation, we consider a fractional derivative of distributed order defined by
\[
\BD_t^{(\mu)}v(t):=\int_0^1\mu(\al)\pa_t^\al v(t)\,\rd\al,
\]
where $\pa_t^\al$ is the Caputo derivative given by
\[
\pa_t^\al v(t)=\f1{\Ga(1-\al)}\int^t_0(t-s)^{-\al}\f{\rd v}{\rd s}(s)\,\rd s.
\]
Throughout this chapter, let $\Om$ be an open bounded domain in $\BR^d$ with smooth boundary $\pa\Om$ of $C^2$-class, and $T>0$ be arbitrarily fixed. We understand that $d\ge1$, that is, we discuss general spatial dimensions, if we do not specify.

Then an equation with time-distributed order derivative is described by
\begin{equation}\label{equ-gov}
\BD_t^{(\mu)}u(x,t)=\tri u(x,t),\quad x\in\Om,\ 0<t<T.
\end{equation}
Here, $u(x,t)$ denotes the density of a substance such as contaminants at the location $x$ and the time $t$.

In order to make \eqref{equ-gov} a model which can interpret observation data better, and well realize asymptotic behavior of the density $u(x,t)$, we have to choose a weight function $\mu(\al)$, $0<\al<1$. This determination problem is inverse problems, and we are requested to determine a function $\mu(\al)$ or related parameters, and this is the subject of the current chapter.

We mainly discuss an initial-boundary value problem:
\begin{equation}\label{equ-u-distri}
\begin{cases}
\BD^{(\mu)}_t u(x,t)=\tri u(x,t), & x\in\Om,\ 0<t<T,\\
u(x,0)=a(x), & x\in\Om,\\
u(x,t)=g(x,t), & x\in\pa\Om,\ 0<t<T,
\end{cases}
\end{equation}
where $\mu\in C[0,1]$, $\ge0$, $\not\equiv0$.

As for references on the forward problems \eqref{equ-u-distri}, we are restricted to Kochubei \cite{K08,K09}, Li, Luchko, and Yamamoto \cite{LLY17,LLY-FCAA}, Luchko \cite{L-FCAA}, Meerschaert, Nane, and Vellaisamy \cite{MNV11} and the references therein.

Henceforth, for simplicity, we consider only the Laplacian $\tri$ in the fractional diffusion equation and a Dirichlet boundary condition.

Moreover, we consider also the following two special cases of \eqref{equ-u-distri}.\medskip

{\bf 1: Single time-fractional diffusion equation}

We formally choose $\mu=\de(\,\cdot\,-\al)$ where $0<\al<1$ is fixed and $\de(\,\cdot\,-\al)$ is a Dirac delta function at $\al$.

Then the distributed order fractional diffusion equation is reduced to a single-term fractional diffusion equation:
\begin{equation}\label{equ-u-single}
\begin{cases}
\pa_t^\al u(x,t)=\tri u(x,t), & x\in\Om,\ 0<t<T,\\
u(x,0)=a(x), & x\in\Om,\\
u(x,t)=g(x,t), & x\in\pa\Om,\ 0<t<T.
\end{cases}
\end{equation}

The diffusion equation with a single-term time fractional derivatives has attracted great attention in different areas, and here we refer to a restricted number of works on the direct problem such as the well-posedness of \eqref{equ-u-single} and other qualitative properties, which are necessary for discussions of the inverse problems: Gorenflo, Luchko, and Yamamoto
\cite{GLY15}, Gorenflo, Luchko, and Zabrejko \cite{GLZ99}, Luchko \cite{Lu09}, Sakamoto and Yamamoto \cite{SY11}, Zacher \cite{Zacher}.

A natural extension for the single-term time-fractional diffusion equation is a multi-term time-fractional diffusion equation:\medskip

{\bf 2: Multi-term time-fractional diffusion equation}

We set $\mu=\sum_{j=1}^\ell p_j\de(\,\cdot\,-\al_j)$ where $0<\al_\ell<\cdots<\al_1<1$, $p_1,\ldots,p_\ell$ are positive constants:
\begin{equation}\label{equ-u-multi}
\left\{\!\begin{alignedat}{2}
& \sum_{j=1}^\ell p_j\pa_t^{\al_j}u(x,t)=\tri u(x,t), & \quad & x\in\Om,\ 0<t<T,\\
& u(x,0)=a(x), & \quad & x\in\Om,\\
& u(x,t)=g(x,t), & \quad & x\in\pa\Om,\ 0<t<T.
\end{alignedat}\right.
\end{equation}

As minimum references on the well-posedness for \eqref{equ-u-multi}, we list Li, Liu, and Yamamoto \cite{LLY-AMC} and
Li, Huang, and Yamamoto \cite{LHY18}.\medskip

For nonstationary partial differential equations such as classical diffusion and wave equations, the inverse problems have been widely studied. More precisely assuming, for example, that $p(x)$ is unknown in the diffusion equation $\pa_t u(x,t)=\rdiv(p(x)\nb u(x,t))$, we are requested to determine $p(x)$, $x\in\Om$ by some limited extra data on the boundary (e.g., Beilina and Klibanov \cite{BK}, Bellassoued and Yamamoto \cite{BY}, Isakov \cite{Is}, Klibanov and Timonov \cite{KT}). See the chapter ``Inverse problems of determining coefficients of the fractional partial differential equation'' of this handbook as for such inverse coefficient problems.

In this chapter, not for classical partial differential equations, we survey other important types of inverse problems of determining $\mu(\al)$ or related parameters which are considered to govern the anomaly of diffusion and describe an index of the inhomogeneity in heterogeneous media. It is not until such parameters are reasonably determined that we can start discussions on the fundamental issues like the existence of solutions to the initial-boundary value problem \eqref{equ-u-distri}. In practice, one can usually identify or estimate such anomaly indices empirically. On the other hand, mathematical discussions for related inverse problems should be not only interesting but also helpful for determination of parameters by experiments. One of such inverse problems is the determination of the unknown parameters in order to match available data such as $u(x_0,t)$, $0<t<T$ at a monitoring point $x_0\in\Om$.

We note that our survey is far from the perfect, because researchers on inverse problems for fractional partial differential equations are very rapidly developed.

\section{Preliminaries}

In this section, we will set up notations and terminologies, review some of standard facts on the fractional calculus, and introduce several important results related to the forward problem for the time-fractional diffusion equations \eqref{equ-gov}, which are the starting points for further researches concerning the theory of inverse problems.

We recall the Mittag-Leffler function $E_{\al,\be}(z)$ defined by
\[
E_{\al,\be}(z)=\sum_{k=0}^\infty\f{z^k}{\Ga(k\al+\be)},\quad z\in\BC,\quad\al,\be>0
\]
and the asymptotic behavior.

\begin{lemma}\label{lem-ML}
Let $\al\in(0,2)$, $\be>0$, and $\mu\in(\f{\al\pi}2,\min\{\al\pi,\pi\})$. Then for $N\in\BN$,
\[
E_{\al,\be}(z)=-\sum_{k=1}^N\f{z^{-k}}{\Ga(\be-\al k)}+O\left(\f1{|z|^{N+1}}\right)\quad\mbox{with }|z|\to\infty,\ \mu\le|\arg z|\le\pi.
\]
Moreover, there exists a positive constant $C=C(\mu,\al,\be)$ such that
\[
|E_{\al,\be}(z)|\le\f C{1+|z|}\quad\mbox{ with }\mu\le|\arg z|\le\pi.
\]
\end{lemma}

Henceforth, $L^2(\Om)$, $H^1_0(\Om)$, $H^2(\Om)$ are the usual $L^2$-space and the Sobolev spaces of real-valued functions: $L^2(\Om)=\{f;\ \int_\Om|f(x)|^2\,\rd x<\infty\}$, and let $(\,\cdot\,,\,\cdot\,)$ denote the scalar product: $(f,g)=\int_\Om f(x)g(x)\,\rd x$ for $f,g\in L^2(\Om)$. Moreover, we define the Laplacian $-\tri$ with the domain $\cD(-\tri)=H^2(\Om)\cap H^1_0(\Om)$ and an eigensystem $\{(\la_n,\vp_n)\}_{n=1}^\infty$ of $-\tri$ is defined by $\{\vp_n\}_{n\in\BN}\subset\cD(-\tri)$ and $0<\la_1<\la_2\le\ldots$ such that $-\tri\vp_n=\la_n\vp_n$ and $\{\vp_n\}_{n\in\BN}$ is a complete orthonormal system of $L^2(\Om)$.

For $\ga\ge0$, the fractional Laplacian $(-\tri)^\ga$ is defined as
\begin{align*}
\cD((-\tri)^\ga) & :=\left\{f\in L^2(\Om);\ \sum_{n=1}^\infty\left|\la_n^\ga(f,\vp_n)\right|^2<\infty\right\},\\
(-\tri)^\ga f & :=\sum_{n=1}^\infty\la_n^\ga(f,\vp_n)\vp_n.
\end{align*}

\section{Single time-fractional diffusion equation}

In this section, we are concerned with inversion for order in the fractional diffusion equation \eqref{equ-u-single} with homogeneous Dirichlet boundary condition. Let us start with some observations mainly about the asymptotic behavior of the solution as $t\to0$ and $t\to\infty$. By the eigenfunction expansion of $u(x,t)$ in terms of the Mittag-Leffler function $E_{\al,1}$, it was shown in \cite{SY11} that the decay rate of the solutions $u$ to \eqref{equ-u-single} is dominated by $t^{-\al}$ as $t\to\infty$, while $\|u(\,\cdot\,,t)\|_{H^2(\Om)}$ is dominated by $t^{-\al}\|a\|_{L^2(\Om)}$ as $t\to0$. It turns out that the short (or long) asymptotic behavior heavily relies on the fractional orders of the derivatives. On the basis of the asymptotics of the solution as $t\to0$ or $t\to\infty$, it is expected to show formulae of reconstructing the order of fractional derivative in time in the fractional diffusion equation by time history at one fixed spatial point.

Hatano, Nakagawa, Wang and Yamamoto \cite{HNWY} proved

\begin{theorem}[Hatano et al. \cite{HNWY}]\label{thm-HNWY}
{\rm(i)} Assuming that
\[
a\in C_0^\infty(\Om),\quad\tri a(x_0)\ne0.
\]
Then
\[
\al=\lim_{t\to0}\f{t\pa_t u(x_0,t)}{u(x_0,t)-a(x_0)}.
\]

{\rm(ii)} Assuming that
\[
a\in C_0^\infty(\Om),\quad a\ge0\mbox{ or }\le0\mbox{ on }\ov\Om.
\]
Then
\[
\al=-\lim_{t\to\infty}\f{t\pa_t u(x_0,t)}{u(x_0,t)}.
\]
\end{theorem}

In addition, let us mention the work from Meerschaert, Benson, Scheffler, and Baeumer \cite{MBS} in which the space-time fractional diffusion equation
\[
\pa_t^\al u=-(-\tri)^{\ga/2}u,
\]
with $0<\al<1$ and $0<\ga<2$ was considered. As is known, the fractional derivative in time here is usually used to describe particle trapping phenomena, while the fractional space derivative is used to model long particle jumps. These two effects combined together produce a concentration profile with a sharper peak, and heavier tails in the anomalous diffusion. We refer to Tatar and Ulusoy \cite{TU13} for the uniqueness for an inverse problem of simultaneously determining the exponents of the fractional time and space derivatives by the interior point observation $u(x_0,\,\cdot\,)$ in $(0,T)$.

\begin{theorem}[Tatar and Ulusoy \cite{TU13}]\label{thm-TU13}
Let $x_0\in(0,1)$ be any fixed point, and let $u[\al,\ga]$ be the weak solution to
\begin{equation}\label{5}
\begin{cases}
\pa_t^\al u=-(-\tri)^{\ga/2}u, & 0<x<1,\ 0<t<T,\\
u(x,0)=a(x), & 0<x<1,\\
u(x,t)=g(x,t), & x=0,1,\ 0<t<T.
\end{cases}
\end{equation}
Let $u[\be,\eta]$ be the weak solution to \eqref{5} where $\al,\ga$ are replaced by $\be\in(0,1)$ and $\eta\in(0,2)$. If $u[\al,\ga](x_0,t)=u[\be,\eta](x_0,t)$ for all $t\in(0,T)$ and $a\in L^2(0,1)$ satisfies
\[
(a,\vp_n)>0\quad\mbox{for all }n\ge1,
\]
then $\al=\be$ and $\ga=\eta$.
\end{theorem}

In the proof of the above theorem, similarly to Hatano et al. \cite{HNWY}, the time-fractional order can be firstly determined by using a long-time asymptotic behavior of the solution to \eqref{5}, which is a direct result from the property of Mittag-Leffler function in Lemma \ref{lem-ML}. The asymptotics of the eigenvalues $\la_n$ of the one-dimensional Laplacian with the zero Dirichlet boundary condition is then used to show the uniqueness inversion for the spatial order. See also Tatar, Tinaztepe, and Ulusoy \cite{TTU16} as for a numerical reconstruction scheme for $\al$ and $\ga$.

\section{Multi-term time-fractional diffusion equation}

In this section, we survey on multi-term time-fractional diffusion equations. In \eqref{equ-u-multi}, we assume the boundary condition $g=0$ and we investigate inverse problems of identifying fractional orders $\al_j$ and coefficients $p_j$.\smallskip

{\bf Inverse Problem.}\ \ Let $x_0\in\Om$ be fixed and let $I\subset(0,T)$ be a non-empty open interval. Determine the number $\ell$ of fractional orders $\al_j$, fractional orders $\{\al_j\}_{j=1}^\ell$ and positive constant coefficients $\{p_j\}_{j=1}^\ell$ of the fractional derivatives by interior measurement $u(x_0,t)$, $t\in I$.

\begin{remark}
We should mention that the number $\ell$ of the fractional derivatives is also unknown in the above inverse problem.
\end{remark}

For the statement of our main results, we introduce some notation. As an admissible set of unknown parameters, we set $\BR_+=\{p>0\}$ and
\begin{equation}\label{eq-admissible}
\cU_{\mathrm{noc}}=\{(\ell,\bm\al,\bm p);\ \ell\in\BN,\ \bm\al=(\al_1,\ldots,\al_\ell)\in(0,1)^\ell,\ 0<\al_\ell<\ldots<\al_1<1,\ \bm p=(p_1,\ldots,p_\ell)\in\BR^\ell_+\}.
\end{equation}

By means of the eigenfunction expansion of $u$, the unique determination of fractional orders is proved by using the $t$-analyticity of the solution and the strong maximum principle for elliptic equations.

\begin{theorem}[Li and Yamamoto \cite{LY-AA}]\label{thm-LY15}
Assuming that
\[
a\ge0\mbox{ in }\Om,\ a\not\equiv0\mbox{ and }a\in H^{2\ga}(\Om)\cap H_0^1(\Om)\mbox{ with some }\ga>\max\left\{\f d2-1,0\right\}.
\]
Let $u[\ell,\bm\al,\bm p],u[m,\bm\be,\bm q]$ be the weak solutions to \eqref{equ-u-multi} with respect to $(\ell,\bm\al,\bm p),(m,\bm\be,\bm q)\in\cU_{\mathrm{noc}}$. Then for any fixed $x_0\in\Om$,
\[
u[\ell,\bm\al,\bm p](x_0,t)=u[m,\bm\be,\bm q](x_0,t),\quad t\in I
\]
implies
\[
\ell=m,\quad\bm\al=\bm\be,\quad\bm p=\bm q.
\]
\end{theorem}

We describe a numerical method for the reconstruction of orders. For simplicity, we fix the number $\ell$ of the fractional derivatives, and denote by $u[\bm\al,\bm p]$ the unique solution to \eqref{equ-u-multi} with $\bm\al=(\al_1,\ldots,\al_\ell)$ and $\bm p=(p_1,\ldots,p_\ell)$. For the numerical treatment, in Li, Zhang, Jia, and Yamamoto \cite{LZJY13}, Sun, Li, and Jia \cite{SLJ17}, the authors reformulated the inverse problem into an optimization problem
\begin{equation}\label{eq-optimization}
\min_{(\bm\al,\bm p)\in\cU_{\mathrm{oc}}}\|u[\bm\al,\bm p](x_0,\,\cdot\,)-h\|_{L^2(0,T)},
\end{equation}
where $\cU_{\mathrm{oc}}$ is the admissible set \eqref{eq-admissible} with fixed $\ell$ and $h\in L^2(0,T)$ stands for available observation data $u(x_0,t)$, $0<t<T$. The feasibility of this optimization is guaranteed by the following result.

\begin{proposition}
For any observation data $h\in L^2(0,T)$, there exists a minimizer to the optimization problem \eqref{eq-optimization}.
\end{proposition}

The existence of a minimizer is easily verified by taking a minimizing sequence and a routine compactness argument. The key is the Lipschitz continuity of the solution with respect to $(\bm\al,\bm p)$, namely,
\[
\|(u[\bm\al,\bm p]-u[\bm\be,\bm q])(x_0,\,\cdot\,)\|_{L^2(0,T)}\le C\sum_{j=1}^\ell(|\al_j-\be_j|+|p_j-q_j|),
\]
where $(\bm\al,\bm p),(\bm\be,\bm q)\in\cU_{\mathrm{oc}}$ satisfy certain additional assumptions.

Needless to say, this methodology also works for the single-term case. However, for $\ell\ge3$ it becomes extremely difficult to reconstruct smaller orders $\al_j$ in $\bm\al$ because they contribute less to the solution. In order to restore the numerical stability, one can attach some regularization terms to \eqref{eq-optimization}, but here we omit the details.

\section{Distributed order fractional diffusion equation}

We focuses on the determination of the weight function $\mu(\al)$ in \eqref{equ-u-distri}, which is important for experimentally evaluating the characteristic of the diffusion process in heterogeneous medium.\smallskip

{\bf Inverse Problem.}\ \ Let $x_0\in\Om$ be fixed and $I\subset(0,T)$ be a non-empty open interval. Let $u$ be the solution to the initial-boundary value problem \eqref{equ-u-distri}. We will investigate whether $u$ in $\{x_0\}\times I$ can determine $\mu$ in $(0,1)$.\smallskip

In this section, we discuss the uniqueness as the fundamental theoretical topic for the inverse problem and attempt to establish results parallel to that for the multi-term case \eqref{equ-u-multi}. We separately discuss two cases: the case of the homogeneous boundary condition and the case of nonhomogeneous boundary condition, and it is technically difficult to consider both simultaneously.

\subsection{Case of the homogeneous boundary condition}

In this subsection, we assume that $g=0$ in \eqref{equ-u-distri}, that is, the boundary value vanishes on $\pa\Om\times(0,T)$. We introduce an admissible set of unknown weight function
\[
\cU_0=\{\mu\in C[0,1];\ \mu\ge0,\not\equiv0\}.
\]

From \cite{LLY17}, the solution to the initial-boundary value problem with the homogeneous boundary condition is $t$-analytic. Therefore, after taking the Laplace transform of the solution $u$ in $t$, we can reformulate the fractional diffusion equation into an elliptic equation with parameters in the Laplace frequency domain. The uniqueness of the inverse problem can be easily derived by noting the strong maximum principle for elliptic equations. We have the following.

\begin{theorem}[Li, Luchko and Yamamoto \cite{LLY17}]\label{thm-LLY17}
We assume that $g=0$ in \eqref{equ-u-distri}. Let $\mu,\om\in\cU_0$ and $u[\mu],u[\om]$ be the solutions to the initial-boundary value problem \eqref{equ-u-distri} with respect to $\mu,\om\in\cU_0$ respectively. Assume that $a\ge0$ in $\Om$, $a\not\equiv0$ and $a\in H^{2\ga}(\Om)\cap H_0^1(\Om)$ with $\ga>\max\{\f d2-1,0\}$. Then $\mu=\om$ if
\[
u[\mu](x_0,t)=u[\om](x_0,t),\quad x_0\in\Om,\ t\in I.
\]
\end{theorem}

The proof of the above theorem heavily relies on the $t$-analyticity of the solution. Let us point out that the proof of the uniqueness for the corresponding inverse problem for the initial-boundary value problem \eqref{equ-u-distri} with nonhomogeneous boundary condition is totally different. In the next subsection, we discuss the case of nonhomogeneous boundary condition.

Finally, as a direct application of well-posedness of the initial-boundary value problem discussed in \cite{LLY17}, we can show the Lipschitz stability of the solutions with respect to weight functions:
\[
\|(u[\mu]-u[\om])(x_0,\,\cdot\,)\|_{L^2(0,T)}\le C\|\mu-\om\|_{L^\infty(0,1)},
\]
where $\mu,\om\in\cU_0$ satisfy certain additional assumptions. We can formulate the minimization problem to prove the existence of a minimizer similarly to Section 4, but we omit the details.

\subsection{Case of nonhomogeneous boundary condition}

In this section, we continue the discussion of the inverse problem of determining the weight function in \eqref{equ-u-distri} but with nonhomogeneous boundary condition.

In the case of $\Om=(0,1)$, we first introduce Rundell and Zhang \cite{RZ16} which investigates an inverse problem of determining the weight function in \eqref{equ-u-distri} in terms of the value of the solution $u$ in the time interval $(0,\infty)$.

\begin{theorem}[Rundell and Zhang \cite{RZ16}]\label{thm-RZ16}
In \eqref{equ-u-distri}, we assume that $\Om=(0,1)$, $a(x)=0$ for $0<x<1$ and $u(1,t)=0$, $u(0,\,\cdot\,)\in L^\infty(0,\infty),\not\equiv0$. Let $x_0\in(0,1)$ be arbitrarily chosen. For arbitrarily fixed $c_0>0$ and $0<\be_0<1$, let $\mu,\om\in C^1[0,1]$ be positive unknown functions and satisfy
\[
\mu(\al),\om(\al)\ge c_0>0\quad\mbox{in }(\be_0,1).
\]
Let $u[\mu]$ and $u[\om]$ be the solutions to the initial-boundary value problem \eqref{equ-u-distri} with respect to $\mu$ and $\om$, respectively, and let the Laplace transform of one of one of $u[\mu](0,\,\cdot\,)$ and $u[\om](0,\,\cdot\,)$ satisfy
\begin{equation}\label{8}
\cL\{u(0,t);s\}\ne0,\quad s\in(0,\infty).
\end{equation}
Then $u[\mu](x_0,t)=u[\om](x_0,t)$ for any $t\in(0,\infty)$ implies $\mu=\om$ on $[0,1]$.
\end{theorem}

\begin{remark}
It is not convenient to use the above result to recover the unknown weight function practically because we have to assume \eqref{8}, which means that one of $\mu$ and $\om$ is known.
\end{remark}

How about using measurement of finite time length? As mentioned above, for the homogeneous boundary case, analytic continuation of the solution ensures the equivalence between the overposed data $u(x_0,t)$, $t\in(0,T)$ and $u(x_0,t)$, $t\in(0,\infty)$, which is the key idea of the proof of Theorem \ref{thm-LY15}. Our question is whether it is possible to prove the uniqueness of the inverse problem of determining the weight function without using the analytic continuation because we cannot rely on the $t$-analytic of the solutions for the nonhomogeneous boundary counterpart.

We have the following lemma.

\begin{lemma}\label{lem-geq0}
Let the weight function $\mu\in C[0,1]$ be non-negative, and not vanish on $[0,1]$. Suppose that the non-negative function $u\in C_0^\infty((0,T_0);H^4(\Om))$ satisfies
\[
\BD^{(\mu)}_t u-\tri u\le,\not\equiv0\quad\mbox{in }\Om\times(0,T),
\]
where $0<T<T_0$. Suppose also that $\Om$ is connected, open and bounded. Then for any $x\in\Om$, there exists $t_x\in(0,T)$ such that $u(x,t_x)>0$.
\end{lemma}

Li, Fujishiro, and Li \cite{LFL18} gave an affirmative answer by setting the boundary value $g$ vanishing near the final time $T$, and restricting the weight function in the following admissible set:
\[
\cU_{\mathrm f}:=\{\mu\in C[0,1];\ \mu\ge,\not\equiv0,\ \mu\mbox{ is finitely oscillatory}\}.
\]
Here, we call that a function $\mu$ is finitely oscillatory if for any $c\in\BR$, the set $\{\al;\ \mu(\al)=c\}$ is finite.

\begin{theorem}[Li, Fujishiro and Li \cite{LFL18}]\label{thm-Li17}
Let $\Om\subset\BR^d$, $d\le3$, be connected, open, and bounded, and let $T>0$, $x_0\in\Om$ be arbitrarily chosen. Let $g\in C_0^\infty((0,T);H^{7/2}(\pa\Om))$ satisfy $g\ge,\not\equiv0$ on $\pa\Om\times(0,T)$. We further suppose that $u[\mu]$ and $u[\om]$ are the solutions to the problem \eqref{equ-u-distri} with respect to $\mu,\om\in\cU_{\mathrm f}$. Then $\mu=\om$ on $[0,1]$ if $u[\mu](x_0,\,\cdot\,)=u[\om](x_0,\,\cdot\,)$ in $(0,T)$.
\end{theorem}

\section{Conclusions and open problems}\label{sec-conc}

In this chapter, we are mainly concerned with the uniqueness results on the determination of a weight function $\mu(\al)$ of distributed order derivatives in $t$ and related parameters of the fractional diffusion equation by one interior point observation data. Our key ideas of proofs are based on the properties such as asymptotic behavior and $t$-analyticity of the solutions of the initial-boundary value problem \eqref{equ-u-distri}.

The uniqueness in determining orders $\al_j$, $j=1,\ldots,\ell$ in \eqref{equ-u-single} or \eqref{equ-u-multi} follows as byproduct from the uniqueness results (e.g., Cheng, Nakagawa, Yamazaki, and Yamamoto \cite{CNYY09}, Li, Imanuvilov, and Yamamoto \cite{LIY}, Kian, Okasanen, Soccorsi, and Yamamoto \cite{KOSY}) for inverse coefficient problems which are surveyed in the chapter ``Inverse problems of determining coefficients of the fractional partial differential equations'' of this handbook. Here as a supplement to this chapter, we give a partial list of numerical researches on the reconstruction for fractional orders and related coefficients for fractional differential equations: Chen, Liu, Jiang, Turner, and Burrage \cite{CLJTB16}, Lukashchuk \cite{L11}, Sun, Li, and Jia \cite{SLJ17}, Tatar, Tinaztepe, and Ulusoy \cite{TTU16}, Yu, Jiang, and Qi \cite{YJQ15}, Yu, Jiang, and Wang \cite{YJW16}, Zhang, Li, Chi, Jia, and Li \cite{ZLCJL12}, Zheng and Wei \cite{ZW11} and the references therein.

On the other hand, the studies on inverse problems of the recovery of the fractional orders or the weight function $\mu(\al)$ in the model \eqref{equ-u-distri}, are far from satisfactory since all the publications either assumed the homogeneous boundary condition (Theorems \ref{thm-HNWY}, \ref{thm-TU13}, \ref{thm-LY15}, and \ref{thm-LLY17}), or studied this inverse problem by the measurement in $t\in(0,\infty)$ (Theorem \ref{thm-RZ16}). Although Theorem \ref{thm-Li17} proved that the weight function $\mu$ in the problem \eqref{equ-u-distri} with nonhomogeneous boundary condition can be uniquely determined from one interior point overposed data but with an additional assumption that the unknown $\mu$ lies in the restrictive admissible set $\cU_{\mathrm f}$. It would be interesting to investigate inverse problems by the value of the solution at a fixed time as the observation data.

Finally, as far as the authors know, the stability of the inverse problem of determining the fractional orders remains open.

\bibliographystyle{unsrt}

\end{document}